\documentclass[12pt]{article}
\usepackage{latexsym, amssymb, amsmath, amscd, amsfonts, epsfig, graphicx, colordvi,verbatim,ifpdf}
\usepackage{amsfonts, amsmath, amssymb}
\usepackage{amssymb,amsfonts,amsmath,latexsym,epsfig,cite, psfrag,eepic,color}
\usepackage{amscd,graphics}
\usepackage{latexsym, amssymb,  amsmath,amscd, amsfonts, epsfig, graphicx, colordvi,amsthm}

\usepackage{graphicx}
\usepackage{color}
\usepackage{ifpdf}
\usepackage{fancybox}

\usepackage{float}

\newtheorem{thm}{Theorem}[section]
\newtheorem{pro}[thm]{Proposition}

\newtheorem{lem}[thm]{Lemma}
\newtheorem{core}[thm]{Corollary}

\def\pf{\noindent{\it Proof.} }
\def\qed{\nopagebreak\hfill{\rule{4pt}{7pt}}
\medbreak}

\newcommand{\ospt}{\mathop{\mathrm{ospt}}\nolimits}

\newcommand{\spt}{\mathop{\mathrm{spt}}\nolimits}

\numberwithin{equation}{section}

\def\qed{\nopagebreak\hfill{\rule{4pt}{7pt}}
\medbreak}

\newlength{\boxedparwidth}
  {\begin{center} \begin{tabular}{|@{\hspace{.315in}}c@{\hspace{.15in}}|}
                  \hline \\ \begin{minipage}[t]{\boxedparwidth}
                  \setlength{\parindent}{.25in}}%
  {\end{minipage} \\ \\ \hline \end{tabular} \end{center}}

\begin{document}

\begin{center}

{\Large \bf Nearly Equal Distributions of

the Rank and the Crank of Partitions}
\end{center}

\vskip 5mm

\begin{center}
{  William Y.C. Chen}$^{1}$,    {  Kathy Q. Ji}$^{2}$
and {Wenston J.T. Zang}$^{3}$ \vskip 2mm

 $^{1,2}$Center for Applied Mathematics\\
Tianjin University,  Tianjin 300072, P. R. China\\[6pt]
    $^{3}$Center for Combinatorics, LPMC\\
   Nankai University, Tianjin 300071, P. R. China\\[6pt]

   \vskip 2mm

    $^1$chenyc@tju.edu.cn, $^2$kathyji@tju.edu.cn,  $^3$wenston@mail.nankai.edu.cn
\end{center}

\vskip 6mm

{\it Dedicated to Professor Krishna Alladi on the Occasion of His Sixtieth Birthday}

\vskip 6mm \noindent {\bf Abstract.}

Let $N(\leq m,n)$ denote the number of partitions of $n$ with rank not greater than $m$,
 and let $M(\leq m,n)$ denote the number of partitions of $n$ with crank not greater than $m$. Bringmann and Mahlburg observed that  $N(\leq m,n)\leq M(\leq m,n)\leq N(\leq m+1,n)$ for   $m<0$ and $1\leq n\leq 100$.     They also pointed out that these inequalities can be restated as the existence of a re-ordering $\tau_n$ on the set of partitions of $n$ such that $|\text{crank}(\lambda)|-|\text{rank}(\tau_n(\lambda))|=0$ or $1$ for all partitions $\lambda$ of $n$, that is, the rank and the crank are
 nearly equal distributions over partitions of $n$. In the study of the
  spt-function, Andrews, Dyson and Rhoades proposed a conjecture on the unimodality of the spt-crank, and they   showed that this conjecture is equivalent to the inequality  $N(\leq m,n)\leq M(\leq m,n)$ for $m<0$ and $n\geq 1$. We proved this conjecture by  combiantorial arguments.
 In this paper, we prove the  inequality $N(\leq m,n)\leq M(\leq m,n)$ for   $m<0$ and $n\geq 1$.
 Furthermore, we define a  re-ordering $\tau_n$ of the partitions $\lambda$ of $n$ and show that this re-ordering  $\tau_n$ leads to the nearly equal distribution of
  the rank and the crank. Using the re-ordering $\tau_n$,  we  give a new combinatorial interpretation of  the
function $\ospt(n)$ defined by Andrews, Chan and Kim, which immediately leads to  an upper bound for  $\ospt(n)$ due to  Chan and Mao.

\section{Introduction}

The objective of this paper is to confirm an observation of Bringmann and Mahlburg \cite{Bringmann-Mahlburg-2009} on the
nearly equal distribution of the rank and crank of partitions. Recall that the rank
of a partition was introduced by Dyson \cite{Dyson-1944} as the largest part of the partition
minus the number of parts. The crank of a partition was defined by Andrews
and Garvan \cite{Andrews-Garvan-1988} as the largest part if the partition contains no ones, and otherwise as the
number of parts larger than the number of ones minus the number of ones.

Let $m$ be an integer. For $n\geq 1$, let $N(m,n)$ denote the number of partitions of $n$ with rank $m$,
and for $n>1$, let $M(m,n)$ denote the number of partitions of $n$ with crank $m$. For $n=1$, set
\[M(0,1)=-1,\,M(1,1)=M(-1,1)=1,\,\]
and for $n=1$ and  $m\neq -1,0,1$, set
\[M(m,1)=0.\]

Define  the rank and the crank cumulation functions by
\begin{equation}\label{def-nleq}
N(\leq m,n)=\sum_{r\leq m} N(r,n),
\end{equation}
and
\begin{equation}\label{def-mleq}
M(\leq m,n)=\sum_{r\leq m} M(r,n).
\end{equation}

Bringmann and Mahlburg \cite{Bringmann-Mahlburg-2009} observed that for $m<0$ and $1\leq n\leq 100$,
\begin{equation}\label{ine-main-res}
N(\leq m,n)\leq M(\leq m,n)\leq N(\leq m+1,n).
\end{equation}
For $m=-1$, an equivalent form of the inequality $N(\leq -1,n)\leq M(\leq -1,n)$ for $n\geq 1$ was conjectured by Kaavya \cite{Kaavya-2011}. Bringmann and Mahlburg \cite{Bringmann-Mahlburg-2009} pointed out that  this observation may also be stated in terms of ordered lists of partitions. More precisely, for $1\leq n\leq 100$, there must be some re-ordering $\tau_n$ of partitions $\lambda$ of $n$ such that
 \begin{equation}\label{main-2-pro}
 |\text{crank}(\lambda)|-|\text{rank}(\tau_n(\lambda))|=0 \ \text{or} \  1.
 \end{equation}
Moreover, they noticed  that  using \eqref{main-2-pro},
 one can deduce the following  inequality on the spt-function $\spt(n)$:
 \begin{equation}\label{equ-up-sptn-1}
 \spt(n)\leq  \sqrt{2n}p(n),
 \end{equation}
 where   $\spt(n)$ is the spt-function  defined by Andrews \cite{Andrews-2008} as  the total number of smallest parts in all partitions of $n$. It should be noted that Chan and Mao \cite{Chan-Mao-2014}   conjectured that for $n\geq 5$,
\begin{equation}
\frac{\sqrt{6n}}{\pi}p(n)\leq spt(n) \leq \sqrt{n}p(n).
\end{equation}

In the study of the spt-crank, Andrews, Dyson and Rhoades \cite{Andrews-Dyson-Rhoades-2013}
conjectured that the sequence $\{N_S(m,n)\}_m$ is unimodal for  $n\geq 1$,  where $N_S(m,n)$ is the number of $S$-partitions of size $n$ with spt-crank $m$. They
showed that this conjecture is equivalent to the inequality
 $N(\leq m,n)\leq M(\leq m,n)$ for $m<0$ and $n\geq 1$.   They  obtained the following asymptotic formula for $N(\leq m,n)-M(\leq m,n)$, which implies that  the inequality  holds for fixed $m$ and sufficiently large $n$.

\begin{thm}[Andrews, Dyson and Rhoades]  For any given $m<0$, we have
\begin{equation}\label{equ-asy-n-m}
M(\leq m,n)-N(\leq m,n)\sim-\frac{(1+2m)\pi^2}{96n}p(n) \quad\text{as}\quad n\rightarrow \infty.
\end{equation}
\end{thm}

We have shown  this inequality holds for all $m<0$ and $n\geq 1$ by constructing an injection  in \cite{Chen-Ji-Zang-2015}. More precisely,

\begin{thm} {\rm\cite{Chen-Ji-Zang-2015}}\label{main-old}
For $m<0$ and $n\geq 1$,
\begin{equation}\label{equ-ine-main-old}
N(\leq m,n)\leq M(\leq m,n).
\end{equation}
\end{thm}

 It turns out that our constructive approach in \cite{Chen-Ji-Zang-2015} also applies to the inequality
\begin{equation}\label{BM-CI}
M(\leq m,n)\leq N(\leq m+1,n),
 \end{equation}
 for $m< 0$ and $n\geq 1$.
 It should be noted that Mao \cite{Mao-2014} obtained the following  asymptotic formula for $N(\leq m+1,n)-M(\leq m,n)$ which implies that the inequality  \eqref{BM-CI} holds for any fixed   $m<0$ and sufficiently large $n$.
 \begin{thm}[Mao] For any given $m<0$, we have
\begin{equation}\label{equ-asy-nk-ine}
N(\leq m+1,n)-M(\leq m,n)\sim\frac{\pi}{4\sqrt{6n}}p(n) \quad \text{as} \quad n\rightarrow\infty.
\end{equation}
 \end{thm}

In this paper, we  show that the inequality \eqref{BM-CI} holds for $m<0$ and $n\geq 1$.

\begin{thm}\label{main-1}
For $m<0$ and $n\geq 1$,
\begin{equation}\label{equ-ine-main}
 M(\leq m,n) \leq N(\leq m+1,n).
\end{equation}
\end{thm}

If we list the set of partitions of $n$ in two ways, one by the ranks, and the
other by cranks, then we are led to a re-ordering
$\tau_n$ of the partitions   of $n$. Using the inequalities
\eqref{ine-main-res} for $m<0$ and $n\geq 1$, we show that the rank and the crank are nearly equidistributed over
  partitions of $n$. Since there may be more than one partition with the same
  rank or crank, the aforementioned listings may not be unique. Nevertheless, this does not
  affect the required property of the re-ordering $\tau_n$.
It should be noted that the above description of $\tau_n$ relies on the
two orders of partitions of $n$, it would be interesting to find a
 definition of $\tau_n$ directly on a partition $\lambda$ of $n$.

\begin{thm}\label{main-2}
For $n\geq 1$, let $\tau_n$ be a re-ordering on the set
of partitions of $n$ as defined above. Then for any partition
$\lambda$ of $n$, we have
\begin{equation}\label{condition}
 \text{crank}(\lambda)-\text{rank}(\tau_n(\lambda))=
 \begin{cases}
  0 & \text{if }  \text{crank}(\lambda)=0,\\[3pt]
 0\ \text{or} \  1 & \text{if }  \text{crank}(\lambda)>0,\\[3pt]
 0\  \text{or} \  -1 & \text{if }\text{crank}(\lambda)<0.
 \end{cases}
 \end{equation}
\end{thm}

Clearly, the above theorem implies   relation \eqref{main-2-pro}.  For example, for $n=4$,  a re-ordering $\tau_4$ is illustrated in Table \ref{table-tau4}.

\begin{table}[h] \label{table-tau4}
\[\begin{tabular}{|c|c|c|c|c|}
  $\lambda$ & $\text{crank}(\lambda)$ & $\tau_4(\lambda)$ & $\text{rank}(\tau_4(\lambda)) $ &$\text{crank}(\lambda)-\text{rank}(\tau_4(\lambda))$  \\  \hline
  $(1,1,1,1)$ & $-4$ & (1,1,1,1) & $-3$ & $-1$ \\[5pt]
    (2,1,1) & $-2$ & (2,1,1) & $-1$ & $-1$ \\[5pt]
  (3,1) & 0 & (2,2) & 0 & 0 \\[5pt]
  (2,2) & 2 & (3,1) & 1 & 1 \\[5pt]
  (4) & 4 & (4) & 3 & 1
\end{tabular}\]
  \caption{The re-ordering $\tau_4$.}
  \end{table}

 We find that  the map $\tau_n$ is related to  the  function $\ospt(n)$ defined by Andrews, Chan and Kim \cite{Andrews-Chan-Kim} as  the difference between the first positive crank moment and the first positive rank moment, namely,
\begin{equation}\label{equ-ospt}
\ospt(n)=\sum_{m\geq 0}mM(m,n)-\sum_{m\geq 0}mN(m,n).
\end{equation}

Andrews, Chan and Kim \cite{Andrews-Chan-Kim} derived the following generating function of $\ospt(n)$.

\begin{thm}[Andrews, Chan and Kim]\label{thm-gef-ospt}
We have
\begin{eqnarray}\label{thm-gef-ospt-eq}
\sum_{n\geq 0}\ospt(n)q^n&=&\frac{1}{(q)_\infty}\sum_{i=0}^\infty \left( \sum_{j=0}^\infty q^{6i^2+8ij+2j^2+7i+5j+2}
(1-q^{4i+2})(1 -q^{4i+2 j+3})\right. \nonumber \\
&&\qquad +\left.\sum_{j=0}^\infty q^{6i^2+8ij+2j^2+5i+3j+1}(1-q^{2i+1})(1- q^{4i+2 j+2})\right).
\end{eqnarray}
\end{thm}

Based on  the above generating function, Andrews, Chan and Kim \cite{Andrews-Chan-Kim} proved the positivity of $\ospt(n)$.

\begin{thm}[Andrews, Chan and Kim]
For $n\geq 1$, $\ospt(n)>0$.
  \end{thm}

They also found a combinatorial interpretation of $\ospt(n)$ in terms of  even strings and odd strings of a partition. The following theorem shows that the function $\ospt(n)$ is related to
the re-ordering $\tau_n$.

\begin{thm}\label{thm-com-int}
For $n>1$, $\ospt(n)$ equals the number of partitions $\lambda$ of $n$  such that $
 \text{crank}(\lambda)-\text{rank}(\tau_n(\lambda))=1$.
\end{thm}

It can be seen that $\tau_n((n))=(n)$ for $n>1$,  since the partition $(n)$ has the largest rank and the largest crank among all partitions of $n$. It follows that $\text{crank}((n))-\text{rank}(\tau_n((n)))=1$ when $n>1$.
  Thus Theorem \ref{thm-com-int} implies that    $\ospt(n)>0$ for $n>1$.

 The following upper bound of $\ospt(n)$ can be derived from Theorem \ref{thm-com-int}.

\begin{core}\label{coro-1}
For $n>1$,
\begin{equation}
\ospt(n)\leq \frac{p(n)}{2}-\frac{M(0,n)}{2}.
\end{equation}
\end{core}

It is easily seen that  $M(0,n)\geq 1$ for $n\geq 3$ since $\text{crank}((n-1,1))=0$. Hence    Corollary \ref{coro-1}  implies the following inequality due to Chan and Mao \cite{Chan-Mao-2014}:
For   $ n\geq 3$,
 \begin{equation}
\ospt(n)< \frac{p(n)}{2}.
\end{equation}

  This paper is organized as follows.   In Section 2, we give a proof of Theorem \ref{main-1} with the aid of the combinatorial construction in \cite{Chen-Ji-Zang-2015}. In Section 3, we demonstrate that Theorem \ref{main-2}   follows from Theorem \ref{main-1}.
   Section 4 provides proofs of  Theorem \ref{thm-com-int} and Corollary \ref{coro-1}. For completeness, we
   include a derivation of  inequality \eqref{equ-up-sptn-1}.

\section{Proof of Theorem \ref{main-1}}

 In this section, we give a proof of Theorem \ref{main-1}. To this end,
 we first reformulate the inequality $M(\leq m,n)\leq N(\leq m+1,n)$ for $m<0$ and $n\geq 1$  in terms of the rank-set. Let $\lambda=(\lambda_1,\lambda_2,\ldots, \lambda_\ell)$ be a partition. Recall that the rank-set  of $\lambda$ introduced by Dyson \cite{Dyson-1989}
 is the infinite sequence
 \[ [- \lambda_1, 1-\lambda_2, \ldots, j-\lambda_{j+1},\ldots,\ell-1-\lambda_{\ell}, \ell, \ell+1,\ldots] .\]

Let $q(m,n)$ denote the number of partitions $\lambda$ of $n$ such that $m$
  appears in the rank-set of $\lambda$.  Dyson \cite{Dyson-1989}  established the following relation:
   For $n\geq 1$,
\begin{equation}\label{dyson-eq}
M(\leq m,n)=q(m,n),
\end{equation}
 see also Berkovich and Garvan \cite[(3.5)]{Berkovich-Garvan-2002}.

Let $p(m,n)$ denote the number of partitions of $n$ with rank larger than or equal to $m$, namely,
\[p(m,n) =\sum_{r\geq m}N(r,n).\]

By establishing the relation
\begin{equation}
M(\leq m,n)-N(\leq m,n)=q(m,n)-p(-m,n),
\end{equation}
for $m<0$ and $n\geq 1$, we see that $M(\leq m,n)\geq N(\leq m,n)$  is equivalent to the inequality $q(m,n)\geq p(-m,n)$. This was justified by a number of injections in \cite{Chen-Ji-Zang-2015}.

Similarly, to prove $N(\leq m+1,n)\geq M(\leq m,n)$   for $m<0$ and $n\geq 1$,
 we need the following relation.

\begin{thm}\label{thm-trans}
For $m<0$ and $n\geq 1$,
\begin{equation}\label{equ-trans-main}
N(\leq m+1,n)-M(\leq m,n)=q(-m-1,n)-p(m+2,n).
\end{equation}
\end{thm}

\pf  Since
\[N(\leq m+1,n)=\sum_{r=-\infty}^{m+1}N(r,n)\]
and
\[p(m+2,n)=\sum_{r=m+2}^\infty N(r,n),\]
we get
\begin{equation}\label{pro-pn-pm+2}
N(\leq m+1,n)=\sum_{r=-\infty}^\infty N(r,n)-p(m+2,n).
\end{equation}
In fact,
\[
\sum_{r=-\infty}^\infty N(r,n)=p(n),
\]
so that \eqref{pro-pn-pm+2} takes the form
\begin{equation}\label{lem-N}
N(\leq m+1,n)=p(n)-p(m+2,n).
\end{equation}
On the other hand, owing to the  symmetry
\[M(m,n)=M(-m,n),\]
   due to Dyson \cite{Dyson-1989}, \eqref{dyson-eq} becomes
\[q(-m-1,n)=\sum_{r=m+1}^\infty M(r,n).
\]
Hence
\begin{equation}\label{pro-m-pn-qm-1}
M(\leq m,n)=\sum_{r=-\infty}^m M(r,n)=\sum_{r=-\infty}^\infty M(r,n)-q(-m-1,n).
\end{equation}
But
\[\sum_{r=-\infty}^\infty M(r,n)=p(n),\]
we arrive at
\begin{equation}\label{lem-M}
 M(\leq m,n)=p(n)-q(-m-1,n).
\end{equation}
Substracting  \eqref{lem-M}  from \eqref{lem-N} gives \eqref{equ-trans-main}. This completes the proof.\qed

In view of Theorem \ref{thm-trans}, we see that  Theorem \ref{main-1} is equivalent to the following assertion.

\begin{thm}\label{thm-main-2}
For $m\geq 0$ and $n\geq 1$,
\begin{equation}\label{ine-trans-2}
q(m,n) \geq p(-m+1,n).
\end{equation}
\end{thm}

Let $P(-m+1,n)$ denote the set of partitions counted by $p(-m+1,n)$, that is, the set of partitions of $n$ with rank  not less than $-m+1$, and let $Q(m,n)$ denote the set of partitions counted by $q(m,n)$, that is, the set of partitions $\lambda$ of $n$ such that $m$ appears in the rank-set of $\lambda$.
Then Theorem \ref{thm-main-2} can be interpreted as the existence of
  an injection $\Theta$ from  the set $P(-m+1,n)$ to the set $Q(m,n)$ for $m\geq 0$ and $n\geq 1$.

In \cite{Chen-Ji-Zang-2015}, we have constructed an injection $\Phi$ from the set $Q(m,n)$ to $P(-m,n)$   for  $m\geq 0$ and $n\geq 1$.  It turns out that the injection $\Theta$ in this paper is less involved  than the injection $\Phi$ in \cite{Chen-Ji-Zang-2015}. More specifically,  to construct the injection
   $\Phi$, the  set $Q(m,n)$ is divided into six disjoint subsets $Q_i(m,n)$ ($1\leq i\leq 6$) and   the set $P(-m,n)$  is divided into eight disjoint subsets $P_{i}(-m,n)$  ($1\leq i\leq 8$). For $m\geq 1$, the  injection $\Phi$
  consists of six injections  $\phi_i$  from the set $Q_i(m,n)$  to the set $P_{i}(-m,n)$, where $1\leq i\leq 6$. When $m=0$, the  injection $\Phi$ requires considerations of more cases.
  For the purpose of this paper,
  the set $P(-m+1,n)$ will be divided into three disjoint subsets $P_{i}(-m+1,n)$ $(1\leq i\leq 3)$
 and the set $Q(m,n)$  will be divided into three disjoint subsets $Q_i(m,n)$  $(1\leq i\leq 3)$. For $m\geq 0$, the injection $\Theta$  consists of three injections  $\theta_1$, $\theta_2$ and $\theta_3$,
  where $\theta_1$ is the identity map, and for $i=2,4$, $\theta_i$
  is an injection from  $P_i(-m+1,n)$  to   $Q_{i}(m,n)$.

  To describe the injection $\Theta$,  we shall   represent the partitions in $Q(m,n)$ and $P(-m+1,n)$  in terms of  $m$-Durfee rectangle symbols. As a generalization of a Durfee symbol defined by Andrews \cite{Andrews-2007},
 an $m$-Durfee rectangle symbol of a partition  is defined in  \cite{Chen-Ji-Zang-2015}. Let
 $\lambda$ be a partition of $n$. The $m$-Durfee rectangle symbol of $\lambda$ is
defined as follows:
\begin{equation}
(\alpha,\beta)_{(m+j)\times j}=\left(\begin{array}{cccc}
\alpha_1,&\alpha_2,&\ldots,&\alpha_s\\[3pt]
\beta_1,&\beta_2,&\ldots,&\beta_t
\end{array} \right)_{(m+j)\times j},
\end{equation}
where $(m+j)\times j$ is the $m$-Durfee rectangle of the Ferrers diagram of $\lambda$ and $\alpha$ consists of  columns to the right of the $m$-Durfee rectangle and   $\beta$ consists of  rows below the $m$-Durfee rectangle, see Figure \ref{fig1}. Clearly, we have
 \[m+j\geq \alpha_1\geq \alpha_2\geq \cdots \geq \alpha_s, \quad j\geq \beta_1 \geq \beta_2 \geq \cdots \geq \beta_t,\]
 and
\[n=j(m+j)+\sum_{i=1}^s \alpha_i+\sum_{i=1}^t \beta_i.\]
 \input{Fig1.TpX}

 When $m=0$, an $m$-Durfee rectangle symbol reduces to a
 Durfee symbol.
 For the partition $\lambda=(7,7,6,4,3,3,2,2,2)$ in Figure \ref{fig1}, the $2$-Durfee rectangle symbol of $\lambda$ is
 \[ \left(\begin{array}{ccccc }
4,&3,&3, &2\\[3pt]
3,&2,&2,&2
\end{array}
\right)_{5\times 3}.\]

Notice that for a partition $\lambda$ with  $\ell(\lambda)\leq m$, where $\ell(\lambda)$ denotes the number of parts of $\lambda$,  it has no $m$-Durfee rectangle. In this case,
 we adopt the convention that the  $m$-Durfee rectangle has
 no columns, that is, $j=0$, and so the $m$-Durfee rectangle symbol of  $\lambda$ is  defined to be $(\lambda',\emptyset)_{m\times  0},$ where $ \lambda'$ is
  the conjugate of $\lambda$. For example, the $3$-Durfee rectangle symbol of $\lambda=(5,5,1)$ is
\[ \left(\begin{array}{ccccc }
3,&2,&2, &2,&2\\[3pt]
&
\end{array}
\right)_{3\times 0}.\]

The partitions in $P(-m+1,n)$ can be characterized in terms of
 $m$-Durfee rectangle symbols.

\begin{pro}\label{pro-dufrec-2} Assume that $m\geq 0$ and $n\geq 1$. Let $\lambda$ be a partition of $n$  and let $(\alpha,\beta)_{(m+j)\times j}$ be the $m$-Durfee rectangle symbol of $\lambda$. Then the rank of $\lambda$ is not less than $-m+1$ if and only if either $j=0$ or  $j\geq 1$ and $\ell(\beta)+1\leq \ell(\alpha).$
  \end{pro}

\pf The proof is substantially the same as that of \cite[Proposition 3.1]{Chen-Ji-Zang-2015}. Assume that the rank of $\lambda$ is not less than $-m+1$. We are going to show that
 either $j=0$ or $j\geq 1$ and $\ell(\beta)+1\leq \ell(\alpha).$ There are two cases:

\noindent Case 1: $\ell(\lambda)\leq m$. We have $j=0$.

\noindent Case 2: $\ell(\lambda)\geq m+1$. We have $j\geq 1$, $\lambda_1=j+\ell(\alpha)$ and $\ell(\lambda)=m+j+\ell(\beta)$. It follows that
\[\lambda_1-\ell(\lambda)=(j+\ell(\alpha))-(j+m+\ell(\beta))=
\ell(\alpha)-\ell(\beta)-m.\]
Since $\lambda_1-\ell(\lambda)\geq -m+1$, we have $\ell(\alpha)-\ell(\beta)\geq 1$, that is,  $\ell(\beta)+1\leq \ell(\alpha).$

Conversely, we assume that
$j=0$ or $j\geq 1$ and $\ell(\beta)+1\leq \ell(\alpha)$.
We proceed to show that the rank of $\lambda$ is not less than $-m+1$.
There are two cases:

\noindent Case 1: $j=0$. Clearly,
  $\ell(\lambda)\leq m$, which implies that  the rank of $\lambda$ is not less than $-m+1$.

\noindent Case 2:  $j\geq 1$ and $\ell(\beta)+1\leq \ell(\alpha)$. Thus we have
 $\lambda_1=j+\ell(\alpha)$ and $\ell(\lambda)=j+m+\ell(\beta)$.
It follows that
\begin{equation}\label{pro3.2-tempt}
\lambda_1-\ell(\lambda)=(j+\ell(\alpha))-(j+m+\ell(\beta))=
-m+(\ell(\alpha)-\ell(\beta)).
\end{equation}
Since $\ell(\alpha)-\ell(\beta)\geq 1$, by \eqref{pro3.2-tempt}, we obtain that
 $\lambda_1-\ell(\lambda)\geq -m+1$. This completes the proof.  \qed

The following  proposition  will be used to describe partitions in $Q(m,n)$ in terms of $m$-Durfee rectangle symbols.

\begin{pro}{\rm \!\!\!\cite[Proposition 3.1]{Chen-Ji-Zang-2015}}\label{pro-dufrec} Assume that
 $m\geq 0$ and $n\geq 1$. Let $\lambda$ be  a partition of $n$ and let $(\alpha,\beta)_{(m+j)\times j}$ be the $m$-Durfee rectangle symbol of $\lambda$. Then $m$ appears in the rank-set of $\lambda$ if and only if either $j=0$ or $j\geq 1$ and $ \beta_1=j.$
  \end{pro}

 If no confusion arises, we do not distinguish the partition $\lambda$ and its $m$-Durfee rectangle
symbol representation. We shall divide the $m$-Durfee rectangle symbols
$(\alpha,\beta)_{(m+j)\times j }$  in
    $P(-m+1,n)$ into three disjoint subsets $P_{1}(-m+1,n)$, $P_{2}(-m+1,n)$ and $P_{3}(-m+1,n)$.
    More precisely,
  \begin{itemize}
\item[(1)] $P_1(-m+1,n)$ is the set of  $m$-Durfee rectangle symbols $(\alpha,\beta)_{(m+j)\times j }$  in $P(-m+1,n)$ for which either of the following conditions holds:\\
     (i) $j=0$;\\
     (ii) $j\geq 1$ and $\beta_1=j$;

\item[(2)] $P_2(-m+1,n)$ is the set of  $m$-Durfee rectangle symbols $(\alpha,\beta)_{(m+j)\times j }$ in $P(-m+1,n)$  such that $j\geq 1$  and $\beta_1=j-1$;

\item[(3)] $P_3(-m+1,n)$ is the set of  $m$-Durfee rectangle symbols $(\alpha,\beta)_{(m+j)\times j }$ in $P(-m+1,n)$ such that $j\geq 2$  and $\beta_1\leq j-2$.

\end{itemize}

The set of $Q(m,n)$ will be divided into the following three   subsets $Q_1(m,n)$, $Q_2(m,n)$ and $Q_3(m,n)$:

\begin{itemize}
\item[(1)] $Q_1(m,n)$ is the set of $m$-Durfee rectangle symbols $(\gamma,\delta)_{(m+j')\times j'}$
  in $Q(m,n)$
  such that either of the following conditions holds:\\
   (i) $j'=0$;\\
   (ii) $j'\geq 1$ and $\ell(\delta)-\ell(\gamma)\leq -1$;

\item[(2)] $Q_2(m,n)$ is the set of $m$-Durfee rectangle symbols $(\gamma,\delta)_{(m+j')\times j'}$ in $Q(m,n)$ such that $j'\geq 1$, $\ell(\delta)-\ell(\gamma)\geq 0$ and $\gamma_1< m+j'$;

\item[(3)] $Q_3(m,n)$ is the set of $m$-Durfee rectangle symbols $(\gamma,\delta)_{(m+j')\times j'}$ in $Q(m,n)$   such that   $j'\geq 1$,   $\ell(\delta)-\ell(\gamma)\geq 0$ and $\gamma_1=m+j'$.

\end{itemize}

We are now ready to define the injections $\theta_i$ from the set $P_i(-m+1,n)$ to the set $Q_{i}(m,n)$, where $1\leq i\leq 3$.
Since   $P_1(-m+1,n)$ coincides with $Q_1(m,n)$, we set
$\theta_1$   to   the identity map. The following
lemma gives an injection $\theta_2$ from $P_2(-m+1,n)$ to $Q_2(m, n)$.

 \begin{lem}\label{phi-2} For $m\geq 0$ and $n>1$,
there is an injection $\theta_2$ from $P_2(-m+1,n)$ to $Q_2(m, n)$.
\end{lem}

\pf  To define the map $\theta_2$, let
 \[\lambda=\left(\begin{array}{c }
\alpha\\[3pt]
\beta
\end{array}
\right)_{(m+j)\times j}=\left(\begin{array}{cccc }
\alpha_1,&\alpha_2,&\ldots, & \alpha_s\\[3pt]
\beta_1,& \beta_2,&\ldots, & \beta_t
\end{array}
\right)_{(m+j)\times j} \]
be an $m$-Durfee rectangle symbol in $P_2(-m+1,n)$.
From the definition of $P_2(-m+1,n)$, we see that  $s-t\geq 1$, $j\geq 1$, $\alpha_1\leq m+j$ and  $\beta_1=j-1$.

Set
   \[\theta_2(\lambda)=\left(\begin{array}{c }
\gamma\\[3pt]
\delta
\end{array}
\right)_{(m+j')\times j'}=\left(\begin{array}{cccccc }
\alpha_1-1,&\alpha_2-1,&\ldots, &\alpha_s-1& \\[3pt]
\beta_1+1,& \beta_2+1,&\ldots, &\beta_t+1,& 1^{s-t}
\end{array}
\right)_{(m+j)\times j}.\]
Clearly, $\theta_2(\lambda)$ is an $m$-Durfee rectangle symbol of $n$. Furthermore,   $j'=j$, $\ell(\delta)-\ell(\gamma)\geq 0$. Since $\alpha_1\leq m+j$,  we see that $\gamma_1=\alpha_1-1\leq m+j-1<m+j'$.
Noting that $\beta_1=j-1$, we get
 $\delta_1=\beta_1+1=j=j'$. Moreover, $\delta_s=1$ since $s-t\geq 1$. This proves that $\theta_2(\lambda)$ is in $ Q_2(m,n).$

To prove that  $\theta_2$ is an injection, define
\[H(m,n)=\{\theta_2(\lambda)\colon \lambda \in P_2(-m+1,n)\}.\]
Let
\[\mu= \left(\begin{array}{cc}
 \gamma\\[3pt]
\delta
\end{array}\right)_{(m+j')\times j'}=\left(\begin{array}{ccccccc}
\gamma_1,&\gamma_2,&\ldots, &  \gamma_{s'}\\[3pt]
\delta_1,& \delta_2,&\ldots, & \delta_{t'}
\end{array}
\right)_{(m+j')\times j'}\]
be an $m$-Durfee rectangle symbol in $H(m,n)$. Since $\mu \in Q_2(m,n)$, we have $t'\geq s'$, $\gamma_1<m+j'$ and $\delta_1=j'$. According to the construction of $\theta_2$,
 $\delta_{t'}=1$. Define
 \[\sigma(\mu)=\left(\begin{array}{c }
\alpha\\[3pt]
\beta
\end{array}
\right)_{(m+j)\times j}=\left(\begin{array}{ccccccc}
\gamma_1+1,&\gamma_2+1,&\ldots, &  \gamma_{s'}+1,& 1^{t'-s'}\\[3pt]
\delta_1-1,& \delta_2-1,&\ldots, & \delta_{t'}-1
\end{array}
\right)_{(m+j')\times j'}.
 \]
Clearly, $\ell(\beta)<t'$ since $\delta_{t'}=1$, so that $\ell(\alpha)-\ell(\beta)\geq 1$. Moreover, since $\delta_1=j'$ and $j'=j$, we see that $\beta_1=\delta_1-1=j'-1=j-1$. It is easily checked that $\sigma(\theta_2(\lambda))=\lambda$ for any $\lambda$ in $P_2(-m+1,n)$. Hence the map $\theta_2$  is an injection from $P_2(-m+1,n)$ to $Q_2(m,n)$. This completes the proof.\qed

For example, for $m=2$ and $n=35$, consider the following $2$-Durfee rectangle symbol in $P_2(-1,35)$:
\[
\lambda=\left(\begin{array}{ccccccc}
5,&5,&3, &  1,& 1\\[3pt]
2,& 2,&1&
\end{array}
\right)_{5\times 3}.
\]
Applying the injection $\theta_2$ to $\lambda$, we obtain
\[
\mu=\theta_2(\lambda)=\left(\begin{array}{ccccccc}
4,&4,&2 &  & \\[3pt]
3,& 3,&2,&1,&1
\end{array}
\right)_{5\times 3},
\]
which is in $Q_2(2,35)$. Applying $\sigma$ to $\mu$, we recover $\lambda$.

The following
lemma gives an injection $\theta_3$ from $P_3(-m+1,n)$ to $Q_3(m, n)$.

 \begin{lem}\label{phi-3} For $m\geq 0$ and $n>1$,
there is an injection $\theta_3$ from $P_3(-m+1,n)$ to $Q_3(m, n)$.
\end{lem}

\pf  Let
 \[\lambda=\left(\begin{array}{c }
\alpha\\[3pt]
\beta
\end{array}
\right)_{(m+j)\times j}=\left(\begin{array}{cccc }
\alpha_1,&\alpha_2,&\ldots, & \alpha_s\\[3pt]
\beta_1,& \beta_2,&\ldots, & \beta_t
\end{array}
\right)_{(m+j)\times j} \]
be an $m$-Durfee rectangle symbol in $P_3(-m+1,n)$.
By definition,  $s-t\geq 1$, $j\geq 2$ and  $\beta_1\leq j-2$.

Define
   \[\theta_3(\lambda)=\left(\begin{array}{c }
\gamma\\[3pt]
\delta
\end{array}
\right)_{(m+j')\times j'}=\left(\begin{array}{cccccccc}
m+j-1,&\alpha_1-1,&\alpha_2-1,&\ldots, &  \alpha_{s}-1&\\[3pt]
j-1,&\beta_1+1,& \beta_2+1,&\ldots, & \beta_{t}+1&,1^{s-t+1}
\end{array}
\right)_{(m+j-1)\times (j-1)}.\]
Evidently, $\ell(\delta)=s+2$ and $\ell(\gamma)\leq s+1$, and so $\ell(\delta)-\ell(\gamma)\geq 1$.
 Moreover, we have $\gamma_1=m+j-1=m+j'$, $\delta_1=j-1=j'$ and
\begin{eqnarray*}
&&\hspace{-0.7cm} j'(m+j')+\sum_{i=1}^{s+1} \gamma_i+\sum_{i=1}^{s+2} \delta_i\\[3pt]
\hspace{-0.5cm}&&=(m+j-1)(j-1)+\left(m+j-1+\sum_{i=1}^{s} (\alpha_i-1)\right)+\left(j-1+ s-t+1+\sum_{i=1}^{t} (\beta_i +1)\right)\\[3pt]
\hspace{-0.5cm}&&=(m+j)j+\sum_{i=1}^{s} \alpha_i+\sum_{i=1}^{t} \beta_i=n.
\end{eqnarray*}
This yields that $\theta_3(\lambda)$ is in $ Q_3(m,n).$ In particular, since $s-t\geq 1$, we see that \begin{equation}\label{equ-delta-s+2-s+1=1}
\delta_{s+2}=\delta_{s+1}=1.
\end{equation}

To prove that the map $\theta_3$ is an injection, define
\[I(m,n)=\{\theta_3(\lambda)\colon \lambda \in P_3(-m+1,n)\}.\]
Let
\[\mu= \left(\begin{array}{cc}
 \gamma\\[3pt]
\delta
\end{array}\right)_{(m+j')\times j'}=\left(\begin{array}{ccccccc}
\gamma_1,&\gamma_2,&\ldots, &  \gamma_{s'}\\[3pt]
\delta_1,& \delta_2,&\ldots, & \delta_{t'}
\end{array}
\right)_{(m+j')\times j'}\]
be an $m$-Durfee rectangle symbol in $I(m,n)$. Since $\mu \in Q_3(m,n)$, we have $t'\geq s'$, $\gamma_1=m+j'$ and $\delta_1=j'$. By the construction of $\theta_3$,   $t'-s'\geq 1$. Define
 \[\pi(\mu)=\left(\begin{array}{c }
\alpha\\[3pt]
\beta
\end{array}
\right)_{(m+j)\times j}=\left(\begin{array}{cccc }
\gamma_2+1,&\ldots ,& \gamma_s'+1, &1^{t'-s'-1}\\[3pt]
  \delta_2-1,&\ldots ,& \delta_{t'}-1&
\end{array}
\right)_{(m+j'+1)\times (j'+1)}.
 \]
It follows from \eqref{equ-delta-s+2-s+1=1}  that $\ell(\beta)\leq t'-3$ and $\ell(\alpha)=t'-2$. Therefore, $\ell(\alpha)\geq \ell(\beta)+1$ and $\beta_1=\delta_2-1\leq j'-1=j-2$,
 so that $\pi(\mu)$ is in $Q_3(m,n)$. Moreover, it can be checked that   $\pi(\theta_3(\lambda))=\lambda$ for any $\lambda$ in $P_3(-m+1,n)$. This proves that the map $\theta_3$  is an injection from $P_3(-m+1,n)$ to $Q_3(m,n)$.\qed

For example, for $m=3$ and $n=63$, consider the following $3$-Durfee rectangle symbol in $P_3(-2,63)$:
\[
\lambda=\left(\begin{array}{ccccccccccc}
7,&7,&4, &  3,& 3,&2,&1\\[3pt]
2,& 2,&2,&1,&1
\end{array}
\right)_{7\times 4}.
\]
Applying the injection $\theta_3$ to $\lambda$, we obtain
\[
\mu=\theta_3(\lambda)=\left(\begin{array}{ccccccccccc}
6,&6,&6,&3, &  2,& 2,&1\\[3pt]
3,&3,& 3,&3,&2,&2,&1,&1,&1
\end{array}
\right)_{6\times 3},
\]
which is in $Q_3(3,63)$. Applying $\pi$ to $\mu$, we recover $\lambda$.

Combining the bijection $\theta_1$ and the injections $\theta_2$ and $\theta_3$,
we are led to an injection $\Theta$ from $P(-m+1,n)$ to $Q(m,n)$, and hence the proof of
Theorem \ref{thm-main-2} is complete. More precisely, for
a partition $\lambda$, define
\[\Theta(\lambda)=\begin{cases}
\theta_1(\lambda), \quad \text{if} \quad \lambda \in P_1(-m+1,n),\\[6pt]
\theta_2(\lambda), \quad \text{if} \quad \lambda \in P_2(-m+1,n),\\[6pt]
\theta_3(\lambda), \quad \text{if} \quad \lambda \in P_3(-m+1,n).
 \end{cases}\]

 \section{Proof of Theorem \ref{main-2}}

 In this section, with the aid of the inequalities \eqref{equ-ine-main-old} in Theorem \ref{main-old} and \eqref{equ-ine-main} in Theorem \ref{main-1}
 for $m<0$ and $n\geq 1$, we show that it is indeed the case that the re-ordering $\tau_n$ leads to
  nearly equal distributions of the rank and the crank.
 For the sake of presentation, the inequalities in  Theorem \ref{main-old} and Theorem \ref{main-1}
 for $m<0$
 can be recast  for $m\geq 0$.

 \begin{thm} For $m\geq 0$ and $n\geq 1$,
\begin{equation}\label{ine-main-res-1}
N(\leq m,n)\geq M(\leq m,n)\geq  N(\leq m-1,n).
\end{equation}
\end{thm}

To see that the inequalities \eqref{ine-main-res-1} for $m\geq 0$ can be derived from  \eqref{equ-ine-main-old}  and  \eqref{equ-ine-main} for $m<0$, we assume that $m\geq 0$, so that
  \eqref{equ-ine-main-old}  and \eqref{equ-ine-main} take the following forms
\begin{equation}
N(\leq -m-1,n)\leq M(\leq -m-1,n)\leq N(\leq -m,n),
\end{equation}
and hence
\begin{equation}\label{ine-main-res-cor-1}
p(n)-N(\leq -m-1,n)\geq p(n)-M(\leq -m-1,n)\geq p(n)-N(\leq -m,n).
\end{equation}
It follows that
\begin{equation}\label{ine-main-res-cor-2}
\sum_{r=-m}^\infty N(r,n)\geq \sum_{r=-m}^\infty M(r,n)\geq \sum_{r=-m+1}^\infty N(r,n).
\end{equation}
Now, by the symmetry $N(m,n)=N(-m,n)$, see   \cite{Dyson-1969}, we have
\begin{equation}\label{ine-main-res-cor-aaa}
\sum_{r=-m}^\infty N(r,n)=N(\leq m,n) \quad \text{and} \quad \sum_{r=-m+1}^\infty N(r,n)= N(\leq m-1,n).
\end{equation}
Similarly,    the symmetry $M(m,n)=M(-m,n)$, see   \cite{Dyson-1989},  leads to
\begin{equation}\label{ine-main-res-cor-bbb}
  \sum_{r=-m}^\infty M(r,n)=M(\leq m,n).
\end{equation}
Substituting \eqref{ine-main-res-cor-aaa} and \eqref{ine-main-res-cor-bbb} into \eqref{ine-main-res-cor-2}, we obtain \eqref{ine-main-res-1}. Conversely, one can reverse the above steps to
derive \eqref{equ-ine-main-old}  and  \eqref{equ-ine-main} for $m<0$ from  \eqref{ine-main-res-1} for $m\geq 0$. This yields that the inequalities \eqref{ine-main-res-1} for $m\geq 0$ are equivalent to the inequalities \eqref{equ-ine-main-old}  and  \eqref{equ-ine-main} for $m<0$.

We are now ready to prove Theorem \ref{main-2}.

\noindent{\it Proof of Theorem \ref{main-2}.}
 Let $\lambda$ be a partition of $n$, and let $\tau_n(\lambda)=\mu$. Suppose  that $\lambda$ is the $i$-th
 partition of $n$ when the partitions of $n$ are listed in the increasing order of their cranks.
  Meanwhile, $\mu$ is also the $i$-th partition
   in the list of  partitions of $n$ in the increasing order of ranks.    Let $\text{crank}(\lambda)=a$ and $\text{rank}(\mu)=b$, so that
 \begin{equation}\label{ine-m-le-n-i}
 M(\leq a, n)\geq i> M(\leq a-1, n),
 \end{equation}
 and
  \begin{equation}\label{ine-n-le-n-i}
N(\leq b, n)\geq i> N(\leq b-1, n).
 \end{equation}

We now consider  three cases:

\noindent Case 1:    $a=0$. We aim to show that $b=0$. Assume to the contrary  that  $b\neq 0$.  There are two subcases:

 Subcase 1.1: $b<0$. From \eqref{ine-m-le-n-i} and \eqref{ine-n-le-n-i},  we have
$$N(\leq -1,n)\geq N(\leq b,n)\geq i>M(\leq -1,n),$$
 which contradicts the inequality $N(\leq m,n)\leq M(\leq m,n)$ in \eqref{ine-main-res}  with $m=-1$.

  Subcase 1.2: $b>0$. From \eqref{ine-m-le-n-i} and \eqref{ine-n-le-n-i},   we see that
$$M(\leq 0,n)\geq i>N(\leq b-1,n)\geq N(\leq 0,n),$$
which contradicts the inequality $M(\leq m,n)\leq N(\leq m,n)$ in \eqref{ine-main-res-1} with $m=0$.
This completes the proof of Case 1.

\noindent Case 2:   $a<0$. We proceed to show that
 $ b=a $ or $a+1.$ By \eqref{ine-m-le-n-i} and the inequality $M(\leq m,n) \leq N(\leq m+1,n)$ in \eqref{ine-main-res} with $m=a$, we see that
\begin{equation}\label{ine-m-le-n-i=temp-1}
N(\leq a+1,n)\geq  i.
\end{equation}
Combining \eqref{ine-n-le-n-i} and \eqref{ine-m-le-n-i=temp-1}, we deduce that
\[N(\leq a+1,n)>N(\leq b-1,n),\]
  and thus
  \begin{equation}\label{ine-re-a+1-b}
  a+1\geq b.
 \end{equation}
 On the other hand, by \eqref{ine-m-le-n-i} and the inequality $N(\leq m,n) \leq M(\leq m,n)$ in \eqref{ine-main-res} with $m=a-1$, we find that
  \[ N(\leq a-1,n)< i.\]
Together with \eqref{ine-n-le-n-i}, this gives
  \[N(\leq a-1,n)< N(\leq b,n),\]
so that  $ a\leq b$. In view of \eqref{ine-re-a+1-b}, we obtain that $b=a$ or $a+1$. This completes the proof
of Case 2.

\noindent Case 3: $a>0$. We claim that
$b=a$   or $a-1.$ Combining the inequality $M(\leq m,n)\geq N(\leq m-1,n)$ in  \eqref{ine-main-res-1}
with $m=a-1$ and  the inequality $M(\leq a-1, n)<i$ in \eqref{ine-m-le-n-i},  we get
\begin{equation}\label{ine-sum-r=a-1-2}
 N(\leq a-2,n)<i.
\end{equation}
By means of  \eqref{ine-n-le-n-i} and \eqref{ine-sum-r=a-1-2}, we find that
\begin{equation*}\label{ine-sum-r=a+1-1}
N(\leq b,n)>N(\leq a-2,n),
\end{equation*}
whence
\begin{equation}\label{ine-a-1-leq-b}
a-1\leq b.
\end{equation}

On the other hand, combining  the inequality $N(\leq m,n)\geq M(\leq m,n)$ in  \eqref{ine-main-res-1} with $m=a$  and the inequality $M(\leq a-1, n)<i$ in  \eqref{ine-m-le-n-i},  we are led to
\begin{equation}\label{ine-sum-r=a-1}
N(\leq a,n)\geq i,
\end{equation}
which together with   \eqref{ine-n-le-n-i}  yields that
\begin{equation*}\label{ine-sum-r=a+1-1}
N(\leq a,n)>N(\leq b-1,n),
\end{equation*}
and hence $a\geq b.$
But it has been shown that $b\geq a-1$, whence the conclusion that $b=a-1$ or $a$. \qed

\section{Proofs of Theorem  \ref{thm-com-int} and Corollary \ref{coro-1} }

In this section, we give a proof of Theorem \ref{thm-com-int} on the
 interpretation of the $\ospt$-function.
      Then we use Theorem \ref{thm-com-int} to deduce Corollary \ref{coro-1},
      which gives an upper bound of the $\ospt$-function. Finally, for
      completeness,
      we include a derivation of \eqref{equ-up-sptn-1} from \eqref{main-2-pro}, as suggested by Bringmann and Mahlburg.

{\noindent \it Proof of Theorem \ref{thm-com-int}.} Let $\mathcal{P}(n)$ denote the set of partitions of $n$.  By the definition \eqref{equ-ospt} of $\ospt(n)$, we see that
\begin{eqnarray}\label{equ-2ospt-1}
\ospt(n)&=&\sum_{\lambda \in \mathcal{P}(n) \atop \text{crank}(\lambda)> 0}\text{crank}(\lambda)-
\sum_{\lambda \in \mathcal{P}(n) \atop \text{rank}(\lambda)>0}\text{rank}(\lambda).
\end{eqnarray}
We claim that
\begin{equation}\label{equ-claim-proof-thm-ospt}
\sum_{\lambda \in \mathcal{P}(n) \atop \text{rank}(\lambda)> 0}\text{rank}(\lambda)=\sum_{\lambda \in \mathcal{P}(n) \atop \text{crank}(\lambda)> 0}\text{rank}(\tau_n(\lambda)).
\end{equation}
From Theorem \ref{main-2},  we see that if $\text{crank}(\lambda)>0,$ then $\text{rank}(\tau_n(\lambda))\geq 0$. This implies that
\begin{eqnarray}\label{equ-inc-lu-sion}
  \{\lambda \in \mathcal{P}(n) \colon \text{crank}(\lambda)>0\}\subseteq\{\lambda  \in  \mathcal{P}(n) \colon \text{rank}(\tau_n(\lambda))\geq 0\}.
\end{eqnarray}
Therefore,
\begin{equation}\label{eqnnnar-aa}
\sum_{\lambda \in \mathcal{P}(n)\atop \text{crank}(\lambda)>0}\text{rank}(\tau_n(\lambda))\leq \sum_{\lambda \in \mathcal{P}(n)\atop \text{rank}(\tau_n(\lambda))\geq 0}\text{rank}(\tau_n(\lambda)).
\end{equation}
From  Theorem \ref{main-2}, we also
 see that if   $\text{crank}(\lambda)=0,$ then $\text{rank}(\tau_n(\lambda))=0$, and if  $\text{crank}(\lambda)<0,$ then $\text{rank}(\tau_n(\lambda))\leq 0$. Now,
\begin{eqnarray}\label{equ-inc-lusion}
\{\lambda \in  \mathcal{P}(n)\colon \text{rank}(\tau_n(\lambda))> 0\}\subseteq \{\lambda \in \mathcal{P}(n) \colon \text{crank}(\lambda)>0\}.
\end{eqnarray}
Hence by \eqref{equ-inc-lu-sion},
\begin{equation}\label{eqnnnar}
\sum_{\lambda \in \mathcal{P}(n)\atop \text{rank}(\tau_n(\lambda))> 0}\text{rank}(\tau_n(\lambda))\leq \sum_{\lambda \in \mathcal{P}(n)\atop \text{crank}(\lambda)>0}\text{rank}(\tau_n(\lambda)).
\end{equation}
Since
\begin{equation*}\label{equ-pf-thm-l-2}
\sum_{\lambda \in \mathcal{P}(n)\atop\text{rank}(\tau_n(\lambda))\geq 0}\text{rank}(\tau_n(\lambda))=\sum_{\lambda \in \mathcal{P}(n)\atop\text{rank}(\tau_n(\lambda))> 0}\text{rank}(\tau_n(\lambda)),
\end{equation*}
from \eqref{eqnnnar-aa} and \eqref{eqnnnar}, we infer that
\begin{equation}\label{equ-pf-con-l-1}
\sum_{\lambda \in \mathcal{P}(n)\atop\text{rank}(\tau_n(\lambda))> 0}\text{rank}(\tau_n(\lambda))=\sum_{\lambda \in \mathcal{P}(n)\atop\text{crank}(\lambda)>0}
\text{rank}(\tau_n(\lambda)).
\end{equation}
But,
\begin{equation}\label{equ-pf-con-l-2}
\sum_{\lambda \in \mathcal{P}(n)\atop\text{rank}(\tau_n(\lambda))> 0}\text{rank}(\tau_n(\lambda))=\sum_{\lambda \in \mathcal{P}(n)\atop\text{rank}(\lambda)> 0}\text{rank}(\lambda).
\end{equation}
Thus we arrive at \eqref{equ-claim-proof-thm-ospt}, and so  the claim is justified.

Substituting \eqref{equ-claim-proof-thm-ospt} into \eqref{equ-2ospt-1}, we get
\begin{eqnarray}\label{equ-2ospt-1-imme}
\ospt(n)
&=&\sum_{\lambda \in \mathcal{P}(n)\atop \text{crank}(\lambda)>0}\text{crank}(\lambda)-
\sum_{\lambda \in \mathcal{P}(n)\atop \text{crank}(\lambda)>0}\text{rank}(\tau_n(\lambda))\nonumber\\[3pt]
&=&\sum_{\lambda \in \mathcal{P}(n)\atop \text{crank}(\lambda)> 0}(\text{crank}(\lambda)-
\text{rank}(\tau_n(\lambda))).
\end{eqnarray}
Appealing to  Theorem \ref{main-2}, we see that if $\text{crank}(\lambda)>0$, then
\[\text{crank}(\lambda)
-\text{rank}(\tau_n(\lambda))=0 \quad \text{or} \quad 1.\]
By \eqref{equ-2ospt-1-imme},
\begin{equation}\label{equ-ospt-co-1-1-1-rr}
\ospt(n)=\#\{\lambda \in \mathcal{P}(n) \colon\text{crank}(\lambda)>0   \text{ and }  \text{crank}(\lambda)
-\text{rank}(\tau_n(\lambda))=1\}.
\end{equation}
Also, by Theorem \ref{main-2}, we see that  if  $\text{crank}(\lambda)
-\text{rank}(\tau_n(\lambda))=1$, then $\text{crank}(\lambda)>0$.
Consequently,
\begin{equation}\label{equ-ospt-co-1-1-1}
\ospt(n)=\#\{\lambda \in \mathcal{P}(n) \colon\text{crank}(\lambda)
-\text{rank}(\tau_n(\lambda))=1\},
\end{equation}
as desired.  \qed

As an application of Theorem \ref{thm-com-int}, we  give a direct proof of  Corollary \ref{coro-1}.

{\noindent \it Proof of  Corollary \ref{coro-1}.}
From the symmetry $M(m,n)=M(-m,n)$, we see that
\begin{equation}\label{ine-sym-mmn}
p(n)=\sum_{m=-\infty}^\infty M(m,n)=M(0,n)+2\sum_{m\geq 1}M(m,n).
\end{equation}
Hence
\begin{equation}\label{equ-po-crank}
\sum_{m\geq 1}M(m,n)=\frac{p(n)}{2}-\frac{M(0,n)}{2}.
\end{equation}
In virtue of Theorem \ref{main-2},  if  $\text{crank}(\lambda)
-\text{rank}(\tau_n(\lambda))=1$, then $\text{crank}(\lambda)>0$, and hence
\begin{equation}\label{equ-ospt-co-1-1-1-cor}
\#\{\lambda \in \mathcal{P}(n) \colon\text{crank}(\lambda)
-\text{rank}(\tau_n(\lambda))=1\}\leq \#\{\lambda \in \mathcal{P}(n) \colon\text{crank}(\lambda)
>0\} .
\end{equation}
This, combined with Theorem \ref{thm-com-int}, leads to
\begin{equation}\label{equ-ospt-po-crank}
\ospt(n)\leq \#\{\lambda \in \mathcal{P}(n) \colon\text{crank}(\lambda)
>0\} = \sum_{m\geq 1}M(m,n).
\end{equation}
Substituting \eqref{equ-po-crank} into \eqref{equ-ospt-po-crank}, we obtain that
\[\ospt(n)\leq \frac{p(n)}{2}-\frac{M(0,n)}{2},\]
as desired.\qed

We conclude  by providing a derivation of  inequality \eqref{equ-up-sptn-1}, that is,
$\spt(n)\leq \sqrt{2n}p(n).$
 Recall that   the $k$-th moment $N_k(n)$  of ranks and the $k$-th  moment   $M_k(n)$   of cranks  were defined by  Atkin and Garvan \cite{Atkin-Garvan-03} as follows:
  \begin{eqnarray}\label{equ-def-nkmn}
N_k(n)&=& \sum_{m=-\infty}^{\infty}m^kN(m,n),\\[3pt]
M_k(n)&=&\sum_{m=-\infty}^\infty m^k M(m,n)\label{equ-def-mk}.
\end{eqnarray}
Andrews \cite{Andrews-2008} showed that the spt-function can be expressed in terms of the second moment $N_2(n)$  of ranks, namely,
\begin{equation}\label{spt-moments}
\spt(n)=np(n)-\frac{1}{2}N_2(n).
\end{equation}
Employing the following relation due to Dyson \cite{Dyson-1989},
\begin{equation}\label{m2-2npn-dyson}
M_2(n)=2np(n),
\end{equation}
Garvan \cite{Garvan-2010} observed that the following expression
\begin{equation}\label{spt-def-eq}
\spt(n)=\frac{1}{2}M_2(n)-\frac{1}{2}N_2(n),
\end{equation}
 implies that $M_2(n)>N_2(n)$ for $n\geq 1$. In general, he
 conjectured and later proved that $M_{2k}(n)>N_{2k}(n)$ for $k\geq 1$ and $n\geq 1$, see \cite{Garvan-2011}.

Bringmann and Mahlburg \cite{Bringmann-Mahlburg-2009} pointed out  that   the   inequality \eqref{equ-up-sptn-1} can be derived by combining  the re-ordering $\tau_n$ and  the Cauchy-Schwarz inequality.
    By \eqref{spt-def-eq}, we see that
\begin{eqnarray}\label{spt-ineq-t--1}
      2\spt(n)&=& \sum_{m=-\infty}^{\infty} m^2 M(m,n)-\sum_{m=-\infty}^{\infty} m^2 N(m,n)\nonumber \\[3pt]
      &=&\sum_{\lambda \in \mathcal{P}(n)}\text{crank}^2(\lambda)-\sum_{\lambda \in \mathcal{P}(n)}\text{rank}^2(\lambda).
      \end{eqnarray}
Since
\[\sum_{\lambda \in \mathcal{P}(n)}\text{rank}^2(\lambda)=\sum_{\lambda \in \mathcal{P}(n)}\text{rank}^2(\tau_n(\lambda)),\]
 \eqref{spt-ineq-t--1} can be rewritten as
\begin{eqnarray}
   2\spt(n)   &=&\sum_{\lambda \in \mathcal{P}(n)}\left(\text{crank}^2(\lambda)-\text{rank}^2(\tau_n(\lambda))\right)
      \nonumber \\[3pt]
      &=&\sum_{\lambda \in \mathcal{P}(n)}(\mid \text{crank}(\lambda)\mid-\mid \text{rank}(\tau_n(\lambda)\mid)\nonumber\\[3pt]
      &&\quad \quad \quad \cdot(\mid \text{crank}(\lambda)\mid+\mid \text{rank}(\tau_n(\lambda))\mid).   \label{spt-ineq-t-1}
\end{eqnarray}
By \eqref{main-2-pro},  we find that
\begin{equation*}\label{ine-c+r-t}
\mid \text{crank}(\lambda)\mid+\mid \text{rank}(\tau_n(\lambda))\mid\leq 2\mid \text{crank}(\lambda)\mid
\end{equation*}
and
\begin{equation*}\label{ine-c+r-t}
0\leq \mid \text{crank}(\lambda)\mid-\mid \text{rank}(\tau_n(\lambda))\mid\leq 1.
\end{equation*}
Thus \eqref{spt-ineq-t-1} gives
\begin{eqnarray}\label{ine-2spt-le-c}
 \spt(n) &\leq& \sum_{\lambda \in \mathcal{P}(n) }\mid \text{crank}(\lambda)\mid.
\end{eqnarray}
Applying the inequality on the arithmetic and quadratic means
\begin{equation}
\frac{x_1+x_2+\cdots+x_n}{n} \leq \sqrt{\frac{x_1^2+x_2^2+\cdots+x_n^2}{n}}
\end{equation}
for nonnegative real numbers to the numbers $ |\text{crank}(\lambda)|$, where
$\lambda$ ranges over partitions of $n$, we are led to
\begin{eqnarray}\label{ine-c-s}
\frac{\sum_{\lambda \in \mathcal{P}(n) }\mid \text{crank}(\lambda)\mid}{p(n)}&&\leq \sqrt{\frac{\sum_{\lambda \in \mathcal{P}(n) }\mid \text{crank}(\lambda)\mid^2}{p(n)}}.\nonumber\\[5pt]
&&=\sqrt{\frac{M_2(n)}{p(n)}}.
\end{eqnarray}
In light of  Dyson's identity \eqref{m2-2npn-dyson}, this  becomes
\begin{eqnarray}\label{ine-c-s-de}
 \sum_{\lambda \in \mathcal{P}(n) }\mid \text{crank}(\lambda)\mid &\leq \sqrt{2n}{p(n)}.
\end{eqnarray}
Combining (\ref{ine-2spt-le-c}) and (\ref{ine-c-s-de}) completes the proof.  \qed

\vskip 0.5cm

\noindent{\bf Acknowledgments.} This work was supported by the 973 Project and the National Science Foundation of China.

\end{document}